\documentclass{article}
\usepackage{amssymb}

\usepackage{graphicx}
\usepackage{amsmath}
\usepackage{a4wide}

\input{tcilatex}

\begin{document}

\title{K-Theory Past and Present}
\author{Michael Atiyah}
\date{(In honour of Friedrich Hirzeburch)}
\maketitle

\thispagestyle{empty}

\section{Early Years}

$K$-theory may roughly be described as the study of additive (or abelian)
invariants of large matrices. \ The key point is that, although matrix
multiplication is not commutative, matrices which act in orthogonal
subspaces do commute. \ Given ``enough room'' we can put matrices $A$ and $B$
into the block form 
\begin{equation*}
\left( 
\begin{array}{ll}
A & 0 \\ 
0 & 1
\end{array}
\right) \quad \left( 
\begin{array}{ll}
1 & 0 \\ 
0 & B
\end{array}
\right)
\end{equation*}
which obviously commute.

Examples of abelian invariants are traces and determinants.

$K$-theory was born during the early days of the famous Bonn Arbeitstagung,
as I shall now briefly recall. \ The prime motivation came from Hirzebruch's
generalisation of the classical Riemann-Roch Theorem. \ This concerns a
complex projective algebraic manifold $X$ and a holomorphic (or algebraic)
vector bundle $E$ over $X.$ \ Then one has the sheaf cohomology groups $%
H^{q}(X,E),$ which are finite-dimensional vector spaces, and the
corresponding Euler characteristic 
\begin{equation*}
\chi (X,E)=\sum_{q=0}^{n}(-1)^{q}\dim H^{q}(X,E)
\end{equation*}
where $n$ is the complex dimension of $X.$ \ One also has topological
invariants of $E$ and of the tangent bundle of $X,$ namely their Chern
classes. \ From these one defines a certain explicit polynomial $T(X,E)$
which by evaluation on $X$ becomes a rational number. \ Hirzebruch's
Riemann-Roch Theorem asserts the equality: 
\begin{equation*}
\chi (X,E)=T(X,E).
\end{equation*}

It is an important fact, easily proved, that both $\chi $ and $T$ are
additive for exact sequences of vector bundles: 
\begin{equation*}
0\rightarrow E^{\prime }\rightarrow E\rightarrow E^{\prime \prime
}\rightarrow 0
\end{equation*}
\begin{eqnarray*}
\chi (E) &=&\chi (E^{\prime })+\chi (E^{\prime \prime }) \\
&& \\
T(E) &=&T(E^{\prime })+T(E^{\prime \prime })
\end{eqnarray*}
where we have dropped the dependence on $X.$

This was the starting point of the Grothendieck generalisation, something
which was expounded at length during the very first Arbeitstagung. \
Grothendieck defined an abelian group $K(X)$ as the universal additive
invariant of exact sequences of algebraic vector bundles over $X,$ so that $%
\chi $ and $T$ both gave homomorphisms of $K(X)$ into the integers (or
rationals).

In fact Grothendieck defined two different $K$-groups, one arising from
vector bundles (denoted by $K^{0})$ and the other using coherent sheaves
(denoted by $K_{0}).$ \ These are formally analogous to cohomology and
homology respectively. \ Thus $K^{0}(X)$ is a ring (under tensor product)
while $K_{0}(X)$ is a $K^{0}(X)$-module. \ Moreover $K^{0}$ is contravariant
while $K_{0}$ is covariant (using a generalisation of $\chi ).$ \ Finally
Grothendieck established the analogue of Poincar\'{e} duality. \ While $%
K^{0}(X)$ and $K_{0}(X)$ can be defined for an arbitrary projective variety $%
X,$ singular or not, the natural map 
\begin{equation*}
K^{0}(X)\rightarrow K_{0}(X)
\end{equation*}
is an isomorphism if $X$ is non-singular.

The Grothendieck Riemann-Roch Theorem concerns a morphism 
\begin{equation*}
f:X\rightarrow Y
\end{equation*}
and compares the direct image of $f$ in $K$-theory and cohomology. \ It
reduces to the Hirzeburch version when $Y$ is a point.

\section{Topological $K$-theory}

Also reported on in the Arbeitstagung was the famous periodicity theorem of
Bott concerning the homotopy of the large unitary groups $U(N)$ (for $%
N\rightarrow \infty ).$

Combining Bott's theorem with the formalism\ of Grothendieck, Hirzebruch and
I, in the late 1950's, developed a $K$-theory based on topological vector
bundles over a compact space. \ Here, in addition to a group $K^{0}(X),$ we
also introduced an odd counterpart $K^{1}(X)$ defined as the group of
homotopy classes of $X$ into $U(N),$ for $N$ large. \ Putting these together 
\begin{equation*}
K^{\ast }(X)=K^{0}(X)\oplus K^{1}(X)
\end{equation*}
we obtained a periodic ``generalised cohomology theory''. \ Over the
rationals the Chern character gave an isomorphism: 
\begin{equation*}
ch:K^{\ast }(X)\otimes Q\cong H^{\ast }(X,Q).
\end{equation*}
But, over the integers, $K$-theory is much more subtle and it has had many
interesting topological applications. \ Most notable was the solution of the
vector fields on spheres problem by Frank Adams, using real $K$-theory
(based on the orthogonal groups).

\section{Relation with Analysis}

Topological $K$-theory turned out to have a very natural link with the
theory of operators in Hilbert space. \ If $H$ is a complex Hilbert space
(of infinite dimension) and $\mathcal{B}(H)$ the space of bounded operators
on $H$ with the uniform norm, then one defines the subspace 
\begin{equation*}
\mathcal{F}(H)\subset \mathcal{B}(H)
\end{equation*}
of\textbf{\ Fredholm operators }$T$ by the requirement that both
Ker\thinspace $T$ and Coker\thinspace $T$ have finite-dimensions. \ The
index is then defined by 
\begin{equation*}
\text{index\thinspace }T=\dim (\text{Ker\thinspace }T)-\dim (\text{%
Coker\thinspace T)}
\end{equation*}
and it has the key property that it is continuous, and therefore constant on
each connected component of $\mathcal{F}(H).$ \ Moreover 
\begin{equation*}
\text{index}:\mathcal{F}\rightarrow Z
\end{equation*}
identifies the components of $\mathcal{F}.$

This has a generalisation to any compact space $X.$ \ To any continuous map $%
X\rightarrow \mathcal{F}$ (i.e. a continuous family of Fredholm operators
parametrized by $X)$ one can assign an index in $K(X).$ \ Moreover one gets
in this way an isomorphism 
\begin{equation*}
\text{index}:[X,\mathcal{F}]\rightarrow K(X)
\end{equation*}
where $[X,\mathcal{F]}$ denotes the set of homotopy classes of maps of $X$
into $\mathcal{F}.$

A\ notable example of a Fredholm operator is an elliptic differential
operator on a compact manifold (these are turned into bounded operators by
using appropriate Sobolev norms).

The Atiyah-Singer index theorem, which computes the index of an elliptic
differential operator (or more generally of a family), is naturally
formulated in terms of $K$-theory and is an extension of the Riemann-Roch
Theorem.

\section{Generalisations}

There are many variants and generalisations of $K$-theory, something which
is not surprising given the universality of linear algebra and matrices. \
In each case there are specific features and techniques relevant to the
particular area.

I will now briefly review some of the many generalisatons. \ First, as
already alluded to, is real $K$-theory based on real vector bundles and the
Bott periodicity theorems for the orthogonal groups: here the period is $8$
rather than $2.$ \ Next there is equivariant theory $K_{G}(X),$ where $G$ is
a a compact Lie group (e.g. a finite group) acting on the space $X.$ \ If $X$
is a point we just get the representation or character ring $R(G)$ of the
group $G.$ \ In general $K_{G}(X)$ \ is a module over $R(G)$ and this can be
exploited in terms of the fixed-point sets in $X$ of elements of $G.$

If we pass from the space $X$ to the ring $C(X)$ of continuous
complex-valued functions on $X$ then $K(X)$ can be defined purely
algebraically in terms of finitely-general projective modules over $X.$ \
This then lends itself to a major generalisation if we replace $C(X),$ which
is a commutative $C^{\ast }$-algebra, by a \textbf{non-commutative} $C^{\ast
}$-algebra. \ This has become a rich theory developed by Kasparov, Connes
and others, which is linked to many basic ideas in functional analysis, in
particular to the von Neumann theory of dimension.

In a different direction, applied to the ring of integers of an algebraic
number field, ideas of Bass, Milnor, Quillen have led to deep connections
with classical number theory.

Finally I should mention the many ways in which $K$-theory has appeared in
physics. \ Part of this is related to the index theorem for Dirac operators
but there are more recent and surprising links to string theory which have
been highlighted by Witten.

Another unexpected link with physics is related to the Verlinde algebra, and
I\ shall now describe this in more detail. \ It is still work in progress
and its full significance has yet to be absorbed.

\section{Projective bundles}

Before getting on to the relation with physics let me first describe an
older generalisation of $K$-theory. \ I start by recalling the connection
between vector bundles and projective bundles. \ Given a vector bundle $V$
over a space $X$ we can form the bundle $P(V)$ whose fibre at $x\in X$ is
the projective space $P(V_{x}).$ \ In terms of groups and principal bundles
this is the passage from $GL(n,C)$ to $PGL(n,C)$ or from $U(n)$ to $PU(n).$
\ We have two exact sequences of groups: 
\begin{equation*}
\begin{array}{lllllllll}
1 & \rightarrow & U(1) & \rightarrow & U(n) & \rightarrow & PU(n) & 
\rightarrow & 1 \\ 
1 & \rightarrow & Z_{n} & \rightarrow & SU(n) & \rightarrow & PU(n) & 
\rightarrow & 1
\end{array}
\end{equation*}
The first gives rise to an obstruction $\alpha \in H^{3}(X,Z)$ to lifting a
projective bundle to a vector bundle, while the second gives an obstruction $%
\beta \in H^{2}(X,Z_{n})$ to lifting a projective bundle to a special
unitary bundle. \ They are related by 
\begin{equation*}
\alpha =\delta (\beta )
\end{equation*}
where 
\begin{equation*}
\delta :H^{2}(X,Z_{n})\rightarrow H^{3}(X,Z)
\end{equation*}
is the coboundary operator. \ This shows that 
\begin{equation*}
n\alpha =0.
\end{equation*}
In fact one can show that any $\alpha \in H^{3}(X,Z)$ of order dividing $n$
arises in this way.

Can we define an appropriate $K$-theory for projective bundles with $\alpha
\neq 0\;?$ \ The answer is yes. \ For each fixed $\alpha $ of finite order
we can define an abelian group $K_{\alpha }(X).$ \ Moreover this is a $K(X)$
module. \ I\ now indicate how these ``twisted'' $K$-groups can be defined.

Note first that, for any vector space $V,$ End\thinspace $V=V\otimes V^{\ast
}$ depends only $P(V).$ \ Hence, given a projective bundle $P$ over $X$ we
can define the associated bundle $\mathcal{E}(P)$ of endomorphism (matrix)
algebras. \ The sections of $\mathcal{E}(P)$ form a non-commutative $C^{\ast
}$-algebra and one can define its $K$-group by using finitely-generated
projective modules. \ This $K$-group turns out to depend not on $P$ but only
on its obstruction class $\alpha \in H^{3}(X,Z)$ and so can be denoted by $%
K_{\alpha }(X).$

In addition to the $K(X)$-module structure of $K_{\alpha }(X)$ there are
multiplications 
\begin{equation*}
K_{\alpha }(X)\otimes K_{\beta }(X)\rightarrow K_{\alpha +\beta }(X).
\end{equation*}

\section{The infinite order case}

The material in \S 5 has been around for many years, but it has only
recently been found useful in physics. \ However in the physical situation
one meets elements $\alpha \in H^{3}(X,Z)$ of infinite order and the
question arises of whether one can still define a ``twisted'' group $%
K_{\alpha }(X).$ \ In fact it is possible to do this and one approach is
being developed by Graeme~Segal and myself. \ I will briefly describe this,
but before I\ begin it is important to emphasise that this infinite order
case brings really new features into play. \ It also involves analysis in an
essential way and is not purely algebraic.

Since an $\alpha $ of order $n$ arises from an obstruction problem involving 
$n$-dimensional vector bundles, it is plausible that, for $\alpha $ of
infinite order, we need to consider bundles of Hilbert spaces. \ But here we
have to be careful not to confuse the ``small'' unitary group 
\begin{equation*}
U(\infty )=\underset{N\rightarrow \infty }{\lim }U(N)
\end{equation*}
with the ``large'' group $U(H)$ of all unitary operators in Hilbert space. \
The small unitary group has interesting homotopy groups given by Bott's
periodicity theorem, but $U(H)$ is contractible (Kuiper's theorem). \ The
means that all Hilbert space bundles (with $U(H))$ as structure group) are
trivial. \ This implies the following homotopy equivalences: 
\begin{eqnarray*}
PU(H) &=&U(H)/U(1)\thicksim CP_{\infty }=K(Z,2) \\
BPU(H) &\thicksim &K(Z,3)
\end{eqnarray*}
where $B$ denotes here the classifying space and on the right we have the
Eilenberg-MacLane spaces. \ It follows that $P(H)$-bundles over $X$ are
classified completely by $H^{3}(X,Z).$ \ Thus, for each \newline
$\alpha \in H^{3}(X,Z),$ there is an essentially unique bundle $P_{\alpha }$
over $X$ with fibre $P(H).$

As in finite dimensions $\mathcal{B}(H)$ depends only on $P(H)$ and so we
can define a bundle $\mathcal{B}_{\alpha }$ of algebras over $X.$ \ We now
let $\mathcal{F}_{\alpha }\subset \mathcal{B}_{\alpha }$ be the
corresponding bundle of Fredholm operators. \ Finally we define 
\begin{equation*}
K_{\alpha }(X)=\text{Homotopy classes of sections of }\mathcal{F}_{\alpha }.
\end{equation*}
This definition works for all $\alpha .$ \ If $\alpha $ is of finite order
then $P_{\alpha }$ contains a finite-dimensional sub-bundle, but if $\alpha $
is of infinite order this is not true. \ Thus we are essentially in an
infinite-dimensional analytic situation.

To get the twisted odd groups we recall that $\mathcal{F}^{1}\subset 
\mathcal{F}$, the space of self-adjoint Fredholm operators, is a classifying
space for $K^{1}$ and so we take $\mathcal{F}_{\alpha }^{1}$ $\subset 
\mathcal{F}_{\alpha }$ define 
\begin{equation*}
K_{\alpha }^{1}(X)=\text{Homotopy classes of sections of }\mathcal{F}%
_{\alpha }^{1}.
\end{equation*}

One peculiar feature of the infinite order case is that all sections of $%
\mathcal{F}_{a}$ lie in the zero-index component, or equivalently that the
restriction map 
\begin{equation*}
K_{\alpha }(X)\rightarrow K_{\alpha }(\text{point)}
\end{equation*}
is zero.

What can we say about the relation between twisted $K$-groups and
cohomology? \ Over the rationals, if $\alpha $ is of finite order, nothing
much changes. \ In particular the Chern character induces an isomorphism. \
However if $\alpha $ is of infinite order something new happens. \ We can
now consider the operation $u\rightarrow \alpha u$ on $H^{\ast }(X,Q)$ as a
differential $d_{\alpha }\;(\alpha ^{2}=0$ since $\alpha $ has odd
dimension). \ We can then form the cohomology with respect to this
differential 
\begin{equation*}
\mathcal{H}_{a}=\text{Ker}\,d_{\alpha }/\func{Im}d_{\alpha }.
\end{equation*}
One can then prove that there is an isomorphism 
\begin{equation*}
K_{\alpha }^{\ast }(X)\otimes Q\cong \mathcal{H}_{\alpha }.
\end{equation*}

\noindent \textbf{Note}: \ In the usual Atiyah-Hirzebruch spectral sequence
relating $K$-theory to integral cohomology all differentials are of finite
order and so vanish over $Q.$ \ In particular $d_{3}=Sq^{3},$ the Steenrod
operation. \ However for $K_{\alpha }$ one finds 
\begin{equation*}
d_{3}u=Sq^{3}u+\alpha u
\end{equation*}
and this explains why an $\alpha $ of infinite order gives the isomorphism
above over the rationals.

Chern classes over the integers are a more delicate matter. \ One can
proceed as follows. \ In $\mathcal{F}$ there are various supspaces $\mathcal{%
F}_{r,s}$ (of finite codimension) where 
\begin{equation*}
\dim \text{Ker}=r\quad \dim \text{Coker}=s
\end{equation*}
and these lie in the component of index $r-s.$ \ Koschorke computed the
cohomology classes $c_{r,s}$ dual to the classes $\mathcal{\bar{F}}_{r,s}$
and identified them with certain determinantal polynomials in the Chern
classes. \ Since the $\mathcal{F}_{r,s}\subset \mathcal{F}$ are invariant
under the action of $U(H)$ it follows that they can be defined fibrewise and
this shows that the classes $c_{r,s}$ can be defined for $K_{\alpha }(X).$ \
However the classes for $r=s$ (and so of index zero) are not sufficient to
generate all Chern classes. \ It is a not unreasonable conjecture that the $%
c_{r,r}$ are the only integral characteristic classes for the twisted $K$%
-theories.

While the use of Hilbert spaces $H$ and the corresponding projective spaces $%
P(H)$ may not come naturally to a topologist, these are perfectly natural in
physics. \ I recall that $P(H)$ is the space of quantum states. \ Bundles of
such arise naturally in quantum field theory.

As I noted earlier Witten has argued convincingly that certain charges in
string-theory are best understood as elements of various $K$-groups, and the
case of the twisted group $K_{\alpha }(X)$ (with $\alpha $ of infinite
order) checks with calculations made in super-gravity.

\section{The Verlinde Algebra}

I shall now describe briefly the very recent work of D. Freed, M. Hopkins
and N. Teleman mentioned earlier.

Let $X=G$ be a compact Lie group and, for simplicity, I shall assume that it
is connected, simple and simply connected, though the theory works in the
general case. We consider $G$ as $G$-space, the group acting on itself 
\textbf{by conjugation.}

Since $H^{3}(G,Z)\cong Z$ we can construct twisted $K$-theories for each
integer $k.$ \ Moreover, we can also do this equivariantly, thus obtaining
abelian groups $K_{G,k}^{\ast }(G).$ \ These will all be $R(G)$-modules.

Now the group multiplication map 
\begin{equation*}
\mu :G\times G\rightarrow G
\end{equation*}
is compatible with conjugation and so is a $G$-map. \ In addition to the
pull back $\mu ^{\ast }$ we can also consider the push-forward $\mu _{\ast
}. $ \ This depends on Poincar\'{e} duality for $K$-theory and it works
also, when appropriately formulated, in the present context.

If $\dim G$ is even this gives us a commutative multiplication on $%
K_{G,k}^{0}(G),$ while for $\dim G$ odd our multiplication is on $%
K_{G,k}^{1}(G).$ \ In either case we get a ring.

The claim of Freed, Hopkins and Teleman is that this ring (according to the
parity of $\dim G)$ is naturally isomorphic to the Verlinde algebra of $G$
at level $k-h$ (where $h$ is the Coxeter number). \ The Verlinde algebra is
a key tool in certain quantum field theories and it has been much studied by
physicists, topologists, group theorists and algebraic geometers. \ The $K$%
-theory approach is totally new and much more direct than most other ways. \
It remains to be thoroughly explored.\vspace{3in}

\noindent Department of Mathematics \& Statistics,

\noindent University of Edinburgh,

\noindent James Clerk Maxwell Buildings,

\noindent King's Buildings,

\noindent Edinburgh \ EH9 3JZ.

\end{document}